\documentclass[12pt,reqno]{amsart}
\usepackage{amsmath, amsthm, amsfonts, amscd,amsfonts,amssymb}

\def\R{\mathbb{R}}

\usepackage[english]{babel}
\usepackage{ifthen}
\usepackage{graphics}
\usepackage{pifont}
\usepackage{float}
\usepackage{array}
\usepackage{hhline}
\usepackage{enumerate}
\usepackage{color}

\def\cyco{\textrm{cyco}\,}

\usepackage{graphicx}
\usepackage{graphics}
\newcommand{\biindice}[3]%
{ #1\begin{array}[t]{c}
{\scriptstyle #2}\\
{\scriptstyle #3}\end{array}
}

\usepackage{pifont}
\usepackage{float}
\usepackage{array}
\usepackage{hhline}
\usepackage{enumerate}
\usepackage{color}

\title{Orthogonal  polynomials and operator orderings}
\author{Adel Hamdi and   Jiang Zeng}
\address{Faculty of Science of Gabes, Department of Mathematics,  Cit\'e Erriadh 6072,  Zrig, Gabes, TUNISIA}
\email{ aadel\_hamdi@yahoo.fr}
\address{Universit\'e de Lyon, Universit\'e Lyon 1, Institute Camille Jordan, UMR 5208 du CNRS, 43, Boulevard du 11 Novembre 1918,
 F-69622 Villeurbanne Cedex, FRANCE}
\email{ zeng@math.univ-lyon1.fr}
\begin{document}
\begin{abstract}
An alternative and combinatorial proof is given for
a connection between a system of Hahn polynomials and identities for symmetric elements in the Heisenberg algebra, which was first observed by Bender, Mead, and Pinsky [Phys. Rev. Lett. 56 (1986),  J. Math. Phys. 28, 509 (1987)] and proved by Koornwinder [J. Phys. Phys. 30(4), 1989].
In the same vein  two results announced by Bender and Dunne [J. Math. Phys. 29 (8), 1988]  connecting a special one-parameter class of Hermitian operator orderings and the continuous Hahn polynomials are also proved.
\end{abstract}
\date{\today}
\maketitle
\section{Introduction}
The Meixner-Pollaczek polynomials are defined by
$$
P_n^{(a)}(x;\phi)=\frac{(2a)_n}{n!}e^{in\phi}{}_2F_1(-n,a+ix;2a;1-e^{-2i\phi}).
$$
Consider the  following special Meixner-Pollaczek polynomial:
\begin{align}
S_n(x)= n!P_n^{(1/2)}\left(\frac{1}{2}x;\frac{1}{2}\pi\right)=i^nn!\sum_{k=0}^{n}\frac{(-1)^{k}}{k!}{n\choose k}\prod_{j=0}^{k-1}(ix+1+2j),
\end{align}
which turns out to be the orthogonal polynomial of degree $n$ on $\R$ with respect to the weight function $x\mapsto 1/\cosh(\pi x/2)$.
Clearly we have
$$
S_{n+1}(x)=xS_{n}(x)-n^2S_{n-1}(x)
$$
with $S_0(x)=1$. The first values of these polynomials are as follows:
\begin{align*}
S_1(x)=x,\quad
S_2(x)=x^2-1,\quad
S_3(x)=x^3-5x,\quad
S_4(x)=x^4-14x^2+9.
\end{align*}
Let $\mathfrak{S}_n$ be  the set of permutations on $\{1,2,\ldots, n\}$.
For any $\sigma\in \mathfrak{S}_n$ let
$\cyco \sigma$ be the number of cycles in
$\sigma$ of odd length. Then it is easy to see that the polynomial $S_n(x)$ has the
following combinatorial interpretation:
$$
S_n(x)=(-i)^n\sum_{\sigma\in \mathfrak{S}_n}(ix)^{\cyco \sigma}.
$$
It is interesting to note that the corresponding moment is the secant number $E_{2n}$ defined by
$$
\sum_{n\geq 0}E_{2n}\frac{x^{2n}}{(2n)!}=\frac{1}{\cos x}.
$$
 If
 \begin{equation}\label{eq:WelAlg}
[q,p]:=qp-pq=i,
\end{equation}
and $T_{n}$ is the sum of all possible terms containing $n$ factors of $p$ and $n$ factors of $q$, then
Bender, Mead and Pinsky~\cite{BMP1, BMP2}  first observed and Koornwinder \cite{Koor}  proved the following result
\begin{equation}\label{eq:koor}
T_{n}=\frac{(2n-1)!!}{n!} S_n(T_{1}).
\end{equation}
For example,  we have $T_{1}=pq+qp$ and $T_{2}=T_{1}^2+p^2q^2+q^2p^2$.
It follows from \eqref{eq:WelAlg} that
$[q,p]^2=(qp)^2-qp^2q-pq^2p+(pq)^2=-1$.
 Hence
$T_{1}^2-1=2((qp)^2+(pq)^2)$ and
\begin{align}\label{eq:calcul}
p^2q^2+q^2p^2+1=p(pq+i)q+q(qp-i)p=(pq)^2+(qp)^2=\frac{1}{2}(T_{1}^2-1).
\end{align}
Finally we get $T_{2}=(T_{1}^2-1)+(p^2q^2+q^2p^2+1)=\frac{3}{2}S_2(T_{1})$.

To prove \eqref{eq:koor},  Koornwinder~\cite{Koor} made use of the connection between Laguerre polynomials and Meixner-Pollaczek polynomials, the Rodrigues formula for Laguerre polynomials, an operational formula involving Meixner-Pollaczek polynomials, and the Schr\"odinger model for the irreducible unitary representations of the three-dimensional Heisenberg group.
 To the best knowledge of the authors Koornwinder's proof is the only published one for \eqref{eq:koor}.
 In this paper we shall give an elementary
 proof of \eqref{eq:koor} using the rook placement interpretation of the normal ordering of two non commutative operators. See 
 \cite{Na}, and  also \cite{Var,BDHPS}  for two recent papers on this theory.

On the other hand,  Bender  and Dunne~\cite{BD} discussed  a correspondence between polynomials and rules for operator orderings.
More precisely, given two operators $q$ and $p$ satisfying \eqref{eq:WelAlg}, they consider the possible operator orderings $O$ as a sum
$$
O(q^np^n)=\frac{\sum_{k=0}^na_k^{(n)}q^kp^nq^{n-k}}{\sum_{k=0}^na_k^{(n)}},
$$
where the coefficients $a_k^{(n)}$ may be chosen arbitrarily with $a_0^{(0)}=a_0^{(1)}=a_1^{(1)}=1$.
Hence $O(q^0p^0)=1$ and $O(q^1p^1)=\frac{1}{2}(qp+pq)$.
They pointed out  that with every operator correspondence rule $O$ one can associate a class of polynomials $p_n(x)$ defined by
\begin{equation}
O(q^np^n)=p_{n}[O(qp)].
\end{equation}

If $a_k^{(n)}={n\choose k}$,  their result is equivalent to
\begin{align}\label{eq:bd1}
O(q^np^n):=\sum_{k=0}^n{n\choose k} q^kp^nq^{n-k}=S_{n}(T_{{1}}).
\end{align}
Their polynomial is actually equal to $2^{-n}S_n(2x)$.
For example, we have
$$
O(q^2p^2)=q^2p^2+2qp^2q+p^2q^2=\frac{T_1^2-3}{2}+2\frac{T_1-i}{2}\cdot \frac{T_1+i}{2}=T_1^2-1.
$$

If
\begin{align}
a_{k}^{(n)}={n+l\choose k}{n+l\choose k+l}{n+l\choose l}^{-1},
\end{align}
where $l$ is an arbitrary parameter, Bender  and Dunne~\cite{BD} observed that the corresponding polynomials belong to the large class of continuous Hahn polynomials.
An explicit formula for the $n$th polynomial is
\begin{align*}
P_{n}(x)=i^n {2n+2l\choose n}^{-1} \frac{\Gamma(n+l+1)}{\Gamma(l+1)} {}_3F_2(-n,n+2l+1,\frac{1}{2}+ix;1,l+1;1).
\end{align*}

Bender and Dunne~\cite{BD}  announced the following result
\begin{align}\label{eq:bd2}
O(q^np^n):= \frac{\sum_{k=0}^n a_{k}^{(n)} q^kp^nq^{n-k}}{\sum_{k=0}^n a_{k}^{(n)}}=
P_{n}(T_{1}/2).
\end{align}
Note that  the denominator has a closed formula
$$
D_n:=\sum_{k=0}^na_k^{(n)}={n+l\choose l}^{-1}\sum_{k=0}^n{n+l\choose k}{n+l\choose n-k}= {n+l\choose l}^{-1}{2n+2l\choose n}.
$$

For example, by \eqref{eq:calcul} we have
\begin{align*}
O(q^2p^2) &= \frac{1}{2}\frac{1+l}{3+2l}(p^2q^2 + 2\frac{2+l}{1+l}qp^2q + q^2p^2)\\
&=\frac{1}{2}\frac{1+l}{3+2l}\left(\frac{T_1^2-3}{2} + \frac{2+l}{1+l}\frac{(T_1+i)(T_1-i)}{2}\right)\\
&=\frac{T_1^2}{4}-\frac{1}{4}\frac{1+2l}{3+2l}=P_2(T_1/2).
\end{align*}
Since \eqref{eq:bd1} and \eqref{eq:bd2} were announced without proof, we shall provide a proof similar to that  of \eqref{eq:koor}.

We shall first recall briefly the rook theory of normal ordering in Section 2 and then prove~\eqref{eq:koor},  \eqref{eq:bd1} and \eqref{eq:bd2} 
in Sections 3, 4 and 5, respectively.

\section{Rook placements and the normal ordering problem}
Let $D$ and $U$ be two operators satisfying  the commutation relation
$[D,U]=1$. Then the algebra generated by $D$ and $U$ is the Weyl algebra. 
Each  element of this algebra, identified as a word $w$ on the alphabet $\{D,U\}$, can be uniquely 
written in the normally ordered  form  as
$$
w=\sum_{r,s}c_{r,s} U^rD^s.
$$
The coefficients $c_{r,s}$ can be computed using the rook theory.  The reader is referred  to 
\cite{Var, BDHPS} for more details.  
Given a word $w$ with $n$ letters $U$'s and $m$ letters $D$'s,
we draw a lattice path from
 $(n, 0)$ to $(0,m)$ as follows:
read the word $\textit{w}$ from  left to right and
draw a unit  line to the right (resp. down) if the letter is $D$ (resp. $U$). This lattice path outlines a Ferrers diagram $B_w$ as 
illustrated  in  Figure~1. 
\begin{center}
\begin{picture}(4,3)(-4,-1.5)
\put(-4,0.5){$w=DUDDUDU\quad \Longrightarrow\qquad $}
\put(0,0){\line(0,1){1.2}}
\put(.4,0){\line(0,1){1.2}}
\put(.8,0){\line(0,1){.8}}
\put(1.2,0){\line(0,2){.8}}
\put(0,0){\line(1,0){1.6}}\put(0,0.4){\line(1,0){1.6}}
\put(0,.8){\line(2,0){1.2}}\put(0,1.2){\line(2,0){.4}}\put(1.6,0){\line(0,2){.4}}
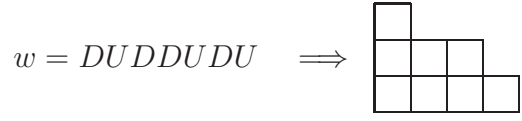
\put(-6,-.8){Figure 1. Correspondence between words and Ferrers diagrams}
\end{picture}
\end{center}

The commutation rule $DU=UD+1$ implies that  the normal writing of 
$w$ amounts  to 
replacing  successively each $DU$ by $UD$ or $1$, this procedure amounts to deleting  each 
{\em up-right-most}  corner  of the Ferrers board    or deleting  it  with its row and column. 
Let  $r_{k}(B)$ be the number of placing $k$ ($k\geq 0$) non-attaking rooks on the Ferrers board $B_w$. It is known (see \cite{Na,Var}) that
\begin{align}\label{eq1}
w=\sum_{k}r_k(B_w) U^{n-k}D^{m-k}.
\end{align}
Now,  it follows from \cite[Theorem 5.1]{Var} that
\begin{align}\label{eq2}
\sum_{B\subseteq [n]\times [n]}r_k(B)=\frac{(2n)!}{2^kk!(n-k)!(n-k)!},
\end{align}
where the sum is over all the Ferrers diagrams contained in the square $[n]\times [n]$ and outlined by the lattice paths 
starting from $(0,n)$ and ending at $(n,0)$.
Let $T_{n}(D,U)$ be the sum of all the words  with $n$ letters $D$ and $n$ letters $U$. We derive from \eqref{eq1} and \eqref{eq2} that
\begin{align}\label{eq:simple2}
T_{n}(D,U)=\sum_{k=0}^n
\frac{(2n)!}{2^kk!(n-k)!(n-k)!}
U^{n-k}D^{n-k}.
\end{align}
On the other hand, let $x:=DU+UD=2UD+1$, then
\begin{align}\label{eq:simple}
UD=\frac{x-1}{2}.
\end{align}
It is also folklore (see \cite[p. 310]{GKP}) that
\begin{align}\label{eq4}
U^nD^n=\sum_{k=0}^n s(n,k) (UD)^k,
\end{align}
where $s(n,k)$ is  the {\em Stirling number of the first kind}.
Since
$$\sum_{j=0}^{n}s(n,j)t^j=t(t-1)\cdots  (t-n+1),
$$
 we can rewrite  \eqref{eq4} as
\begin{align}
U^nD^n=\prod_{j=1}^{n}(UD-j+1).\label{eq:simple1}
\end{align}

\section{Proof of equation~\eqref{eq:koor}}
If we set $D=-iq$ and $U=p$, then equation \eqref{eq:WelAlg} becomes
 $[D,U]=1$ and
 the algebra generated by $D$ and $U$  is
 the Weyl algebra.
Substituting \eqref{eq:simple1} and \eqref{eq:simple} into \eqref{eq:simple2},  we obtain, by replacing $k$ by $n-k$,
\begin{align*}
T_{n}(D,U)&= (2n-1)!!\sum_{k=0}^n   \frac{(-1)^{k}}{k!}
{n\choose k} \prod_{j=0}^{k-1} (-x+2j+1)=(-i)^n\frac{(2n-1)!!}{n!} S_n(ix).
\end{align*}
Since $T_{n}(D,U)=(-i)^nT_{n}(q,p)$, letting $T_1=pq+qp=ix$,  we derive
$$
T_{n}(q,p)=\frac{(2n-1)!! }{n!}S_n(T_1),
$$
which is exactly \eqref{eq:koor}.

\section{Proof of equation~\eqref{eq:bd1}}
For the commutation relation $DU-UD = 1$ and for $n\geq k\geq 0$,  it is readily seen,  by induction on $k$, that
\begin{align}\label{eq:simple4}
D^{k}U^{n} = \sum_{j=0}^k{k\choose j}n^{\underline{j}}U^{n-j}D^{k-j},
\end{align}
where $n^{\underline{j}} = n(n-1)...(n-j+1)=n!/(n-j)!$, then
\begin{align*}
O(D^nU^n)&=\sum_{k=0}^n
{n\choose k} \sum_{j=0}^k{k\choose j}n^{\underline{j}}U^{n-j}D^{n-j}\\
&=\sum_{j=0}^n{n\choose j} 2^{n-j}n^{\underline{j}}
U^{n-j}D^{n-j}.
\end{align*}
Substituting \eqref{eq:simple1} and \eqref{eq:simple} in the last sum,
we obtain, by replacing $n-j$ by $k$,
\begin{align*}
O(D^nU^n)&=\sum_{k=0}^n(-1)^k
{n\choose k} n^{\underline{n-k}}\prod_{j=0}^{k-1} (-x+2j+1)
=\sum_{k=0}^n(-1)^k
{n\choose k} \frac{n!}{k!}\prod_{j=0}^{k-1} (-x+2j+1).
\end{align*}
As $T_{1}= qp + pq = ix$ we derive
$O(q^np^n)=i^nO(D^nU^n)= S_{n}(ix)$.
\section{Proof of equation~\eqref{eq:bd2}}
By definition and \eqref{eq:simple4}, we have
\begin{align*}
O(D^nU^n)
&={2n+2l\choose n}^{-1} \sum_{k=0}^n{n+l\choose k}{n+l\choose k+l}
D^kU^nD^{n-k}\\
&={2n+2l\choose n}^{-1} \sum_{k=0}^n {n+l\choose k}{n+l\choose k+l}\sum_{j=0}^k {k\choose j}n^{\underline{j}}U^{n-j}D^{n-j}\\
&={2n+2l\choose n}^{-1} \sum_{j=0}^n {n\choose j}U^{n-j}D^{n-j} \frac{(n+l)!}{(n+l-j)!}{2n+2l-j\choose n-j}.
\end{align*}
Substituting \eqref{eq:simple1} and
\eqref{eq:simple} and replacing $j$ by $n-k$ we obtain
\begin{align*}
O(D^nU^n)
&={2n+2l\choose n}^{-1} \frac{(n+l)!}{l!} \sum_{k=0}^n {n\choose k} \frac{(n+2l+1)_k}{(l+1)_k} \frac{(-1)^k}{2^kk!} \prod_{j=0}^{k-1}(-x/2+1/2+j).
\end{align*}
Since $O(q^np^n)=i^nO(D^nU^n)$ and $T_{1}=ix$, we derive Eq.~\eqref{eq:bd2}.

\medskip
{\bf Acknowledgments}:
This work was partially
supported by  the French National Research Agency through 
the  grant ANR-08-BLAN-0243-03,  and initiated
 during  the first author's visit to Institut Camille Jordan, Universit\'{e} Lyon 1 in the summer of 2009.
\goodbreak


\end{document}